\documentstyle[amsfonts,12pt]{article}
\newcommand{\Y}{\mbox{$Y$}}
\newcommand{\M}{\mbox{$M$}}

\newcommand{\OM}{\mbox{$\Omega^{1}M$}}

\newcommand{\X}{\mbox{$X$}}

\newcommand{\f}{\mbox{$\cal O$}}
\newcommand{\ko}{{\sf k}}
\newcommand{\codim}{\mbox{codim}}

\newcommand{\blm}{\mbox{$\Lambda$}}

\newcommand{\lon}{\mbox{$\longrightarrow$}}
\newcommand{\ron}{\rightarrow}
\newcommand{\bnu}{\nu}
\newcommand{\bmu}{\mu}

\newcommand{\sym}{{\sf sym}}

\newtheorem{theorem}{Theorem}

\newtheorem{corollary}[theorem]{Corollary}

\begin{document}
\vspace{45 mm}
\title{Geometry of Kodaira moduli spaces }

\author{Sergey A.\ Merkulov\\
\small   School of Mathematics and Statistics, University of Plymouth \\
\small                Plymouth, Devon PL4 8AA, United Kingdom
}

\date{}
\maketitle
\begin{abstract}
A general theorem on the existence of natural torsion-free
affine connections on a complete family of compact complex submanifolds
in a complex manifold is proved.
\end{abstract}

\sloppy

\paragraph{0. Introduction.}
Let $X$ be a compact complex submanifold of a complex manifold $Y$.
The chief object of interest in this paper is the set $M$ of all
holomorphic deformations of $X$ {\em inside} Y, i.e.\ a point $t$
in $M$ can be thought of as a "nearby" compact complex submanifold
$X_t$ in $Y$. The first natural question about this set --- when
is it a manifold? --- was answered by Kodaira \cite{Kodaira1} in 1962.
He found a surprisingly simple sufficient condition, the vanishing
of the first cohomology group $H^1(X,N)$ with coefficients in the normal
sheaf $N\equiv TY\mid_X /TX$, for $M$ to have an induced complex manifold
structure. Moreover, he established a clear correspondence between
the {\em first}\, floors of two towers of infinitesimal neighbourhoods
generated by embeddings,
$t\hookrightarrow M$ and $X_t\hookrightarrow Y$,
by proving that there is a {\em canonical}\, isomorphism
$\ko_{t}: T_{t}\M \longrightarrow  H^{0}(X_{t}, N_{t})$
which associates a global
section of the normal bundle $N_{t}$ of $\X_{t} \hookrightarrow Y$ to any
tangent vector at the corresponding point $t \in \M$. The complete
parameter space $M$ is called a Kodaira moduli space.

 About 14 years later Penrose~\cite{Penrose} studied the Kodaira
moduli space $M$ of rational curves $X_t \simeq {\Bbb C\Bbb P}^{1}$
embedded into a 3-dimensional complex manifold $Y$  with normal bundle
$N_t\cong \Bbb C^{2}\otimes\f(1)$. His striking conclusion was that
4-fold $M$ comes equipped canonically with  an induced half-flat
conformal structure. This Penrose idea to study
geometric structures canonically induced on Kodaira moduli spaces
led to a radically new way of solving of a number non-linear
differential equations in the form of constraints on curvature and
torsion tensors  of an affine connection revealing thereby deep
interconnections between multidimensional complex analysis and differential
geometry (see, e.g., the books \cite{BB,BaE,Besse,Huggett,Manin,PR,WW} and
references cited therein). Instead of analysing some particular Kodaira
moduli spaces (as is normally done in twistor theory, where moduli spaces of
rational curves and quadrics with specific normal bundles have been
only considered), we attempt in this paper
to find sufficient conditions for a general Kodaira moduli space
to have a canonically induced family of torsion-free affine connections
satisfying some natural integrabily conditions. We do it through the analysis
of the {\em second}\, floors of towers of infinitesimal neighbourhoods
generated
by embeddings $t\hookrightarrow M$ and $X_t\hookrightarrow Y$.
Probably the most unexpected aspect of the emerging theory is that both
the Kodaira condition for existence of an induced complex manifold
structure and the conditions for existence of an induced geometric structure
on $M$ have the form of vanishing of some cohomology groups associated
with normal bundles. As an application we show how several "curved twistor"
constructions follow from one and the same general theorem thus merging
several fragments into one picture.

\paragraph{1. Kodaira's existence theorem.}
Let $Y$ and $M$ be complex manifolds and let $\pi_{1}: Y\times  M\lon Y$ and
$\pi_{2} : Y\times M \lon M$ be the natural projections. An analytic family of
compact submanifolds of the complex manifold $Y$ with the moduli space $M$
is a complex submanifold $F\hookrightarrow Y\times M$ such that the restriction
of the projection $\pi_{2}$ on $F$ is a proper regular map (regularity means
that
the rank of the differential of $\bnu\equiv \pi_{2}\mid_{F}: F\lon M$ is equal
at every point to
$\dim M$). Thus the family $F$ has the structure of a double fibration
$$
Y \stackrel{\bmu}{\longleftarrow} F \stackrel{\bnu}{\lon} M
$$
where $\bmu\equiv \pi_{1}\mid_{F}$.
 For each $t\in M$ the compact complex submanifold
$X_{t} = \bmu\circ\bnu^{-1}(t) \hookrightarrow Y$ is said to belong
 to the family $M$.
If $F\hookrightarrow Y\times M$ is an analytic family of compact submanifolds,
then, for any $t\in M$, there is a natural linear map
$$
\ko_{t} : T_{t}M \lon H^{0}(X_{t}, N_t),
$$
 from the tangent space  at $t$ to the vector space of
global holomorphic sections of the normal bundle $N_t =
\left.TY\right|_{X_t}/ TX_{t}$ to the submanifold $X_{t}\hookrightarrow Y$.

An analytic family \{$\X_{t}\hookrightarrow Y \mid t \in \M$\} of compact
submanifolds is
called {\em complete} if the Kodaira map $\ko_{t}$ is an isomorphism at each
point $t$ in the moduli space $M$.
It is called {\em maximal at a point} $t_{0}\in M$,
if, for any other analytic family
 \{$\X_{\tilde{t}}\hookrightarrow Y \mid \tilde{t} \in \tilde{\M}$\}
of compact complex submanifolds such that
$X_{t_{0}} = X_{\tilde{t}_0}$ for a point $\tilde{t}_{0}\in \tilde{M}$,
there exists a neighbourhood $\tilde{U}\subset \tilde{M}$ of
$\tilde{t}_0$ and a holomorphic map
$f: \tilde{U}\lon M$ such that $f(\tilde{t}_{0}) = t_{0}$ and
$X_{f(\tilde{t}')} = X_{\tilde{t}'}$ for each $\tilde{t}'\in \tilde{U}$.
The family  \{$\X_{t}\hookrightarrow Y \mid t \in \M$\} is called {\em maximal}
 if it is maximal at each point $t$ in the moduli space $M$.

In 1962 Kodaira proved that {\em
if $\X \hookrightarrow \Y$ is a compact complex submanifold with normal bundle
$N$ such that $H^{1}(\X,N) = 0$, then \X\ belongs to the
 complete analytic family \mbox{$\{\X_{t} : t \in M\}$} of compact submanifolds
$\X_{t}$ of $Y$. The family is maximal and
its moduli
space is of complex dimension $\dim_{\Bbb C}H^{0}(\X, N)$.
}
\vspace{2 mm}

\paragraph{2. Alpha subspaces of moduli spaces.}
Let \{$\X_{t}\hookrightarrow Y \mid t \in \M$\} be a complete
family of compact complex submanifolds.  For any point
$y\in Y'\equiv \cup_{t\in M}X_t$, there is an associated subset
$\bnu\circ\bmu^{-1}(y)$ in $M$. It is easy to show that
such a subset is always an analytic subspace of $M$. Moreover, if the natural
evaluation map
$$
H^{0}(X_{t}, N_t)\otimes \f_{X_{t}}\lon N_t
$$
is an epimorphism  for all $t\in \bnu\circ\bmu^{-1}(y)$, then
the subspace $\bnu\circ\bmu^{-1}(y)\subset M$ has no singularities, i.e.\
it is a submanifold. In general, we denote the manifold content of
$\bnu\circ\bmu^{-1}(y)$ by $\alpha_{y}$ and, following the terminology of
the seminal paper by Penrose \cite{Penrose}, call it an {\em alpha subspace}
of $M$.

\paragraph{3. Two natural sheaves on moduli space.}
Let $F\hookrightarrow Y\times M$ be a complete analytic family of compact
submanifolds. The Kodaira theorem asserts the existence of a sheaf isomorphism
$\ko : TM\lon \bnu_{*}^0(N_F)$, where $N_F$ denotes the normal bundle of $F$
in $Y\times M$. It is easy to see that $\ko^{-1}$ induces canonically the
following morphisms of $\f_M$-modules
$$
\begin{array}{rccc}
\phi_1: & \bnu_{*}^{0}\left(N_{F}\otimes S^{2} (N_{F}^{*})\right)&
\lon  &TM\otimes S^{2}(\OM), \\
\phi_2: & \bnu_{*}^{0}\left(N_{F}\otimes N_{F}^{*}\right) & \lon &
TM\otimes \OM,
\end{array}
$$
whose images we denote by $\blm$ and $\Theta$ respectively. Thus, the  locally
free sheaves $TM\otimes S^{2}\OM$ and $TM\otimes\OM$ on any Kodaira moduli
space $M$ come equipped canonically with $\f_M$-submodules $\blm$ and,
respectively,  $\Theta$ which may fail to be locally free in general.

\paragraph{4. Vector bundles on $Y$ and affine connections on
$M$.} Let \{$\X_{t}\hookrightarrow Y \mid t \in \M$\} be a complete family
of compact complex submanifolds. After the work by Ward \cite{W1} on instanton
solutions of Yang-Mills equations much attention has been paid  to
holomorphic vector bundles $E$ on $Y$ which are trivial on submanifolds
$X_{t}$ for all $t\in M$. Their geometric
role is that they generate vector bundles on $M$ together with linear
connections which are integrable on alpha subspaces of $M$ (see, e.g.,
\cite{Manin,WW}). In this subsection
we find a geometric meaning of holomorphic vector bundles on $Y$
which, when restricted to $X_{t}$, $t\in M$, are canonically isomorphic to
the normal bundle $N_t$ of the embedding $X_{t}\hookrightarrow Y$, i.e.\
$\left.E\right|_{X_t} = N_t$.

\begin{theorem}\label{2.1}
Let \{$\X_{t}\hookrightarrow Y \mid t \in \M$\} be a complete
family of compact submanifolds.
Suppose there is a holomorphic vector bundle $E$ on $Y$ such that
$\left.E\right|_{X_t} = N_t$ for all $t$ in some domain $M_0\subseteq M$.
If $H^{0}\left(X_t, N_t\otimes S^{2}(N_t^{*})\right) =  H^{1}\left(X_t,
N_t\otimes S^{2}(N_t^{*})\right) = 0$ for each $t\in M_0$, then $M_0$ comes
equipped canonically with an induced torsion-free affine connection such
that, for every $y\in Y'$, the associated
alpha subspace $\alpha_y\cap M_0$ is totally geodesic.
\end{theorem}

 Let $(X, \f_X)$ be an analytic
subspace of a complex manifold $(Y,\f_Y)$ defined by the sheaf of ideals
$J\subset \f_Y$. The {\em m}\,th-order
infinitesimal neighborhood of \X\ in \Y\ is
 the ringed space $\X^{(m)} = (\X, \f^{(m)}_{X})$ with the structure
sheaf $\f^{(m)}_{X} = \f_Y/ J^{m+1}$.
With the $(m+1)$\,th-order infinitesimal neighborhood there is  naturally
associated an {\em mth-order conormal sheaf} of $\f^{(m)}_X$-modules
$\f^{(m)}_{X}(N^{*}) = J/J^{m+2}$ (cf.\ \cite{EL}). By construction,
$\f^{(0)}_{X}(N^{*})$ is the usual conormal sheaf $N^{*}$,
while $\f^{(1)}_{X}(N^{*})$  fits into the exact sequence of
$\f_{X}^{(1)}$-modules
$$
0 \longrightarrow  S^2(N^{*})\longrightarrow \f^{(1)}_{X}(N) \longrightarrow
 N^{*} \longrightarrow 0.
$$
When $X$ is a point in $Y$, then $\left(\f^{(0)}_{X}(N^{*})\right)^{*}$ is
identical to the tangent space $T_X Y$ at $X\in Y$, while
$\left(\f^{(1)}_{X}(N^{*})\right)^{*}$ is usually
denoted by $T_X^{[2]}Y$ and called second-order tangent space.

{\em Proof of Theorem~1}. By assumption,
there is a holomorphic vector bundle on $Y$ such that
$\left. E\right|_{X_t} = N_t$ for any given $t\in M$. The restriction of
its dual to the first order infinitesimal neighbourhood of $X_t$ in $Y$ is an
$\f_{X_t}^{(1)}$-module which fits into the exact sequence \cite{Manin}
$$
0\lon N_t^{*}\otimes \left.E^{*}\right|_{X_t} \stackrel{i}{\lon}
\left. E^{*}\right|_{X_t^{(1)}} \lon
\left.E^{*}\right|_{X_t} \lon 0  .
$$
Therefore, one obtains an extension
$$
0 \longrightarrow  S^2(N_t^{*}) \longrightarrow
\left.E^{*}\right|_{X_t^{(1)}}
/i(\wedge^{2}N^{*})\longrightarrow
N_t^{*} \longrightarrow 0.
$$
It is straightforward to show using local coordinate representations of
all objects involved \cite{MeI}, that the difference
$$
\left[ \f^{(1)}_{X_t}(N^{*}) - \left.E^{*}\right|_{X_t^{(1)}}
/i(\wedge^{2}N^{*})\right]\in \mbox{\rm Ext}_{\f^{(1)}_{X}}\left(N^{*},
S^2(N^{*})\right)
$$
is a locally free $\f_{X_t}$-module which compares to 2nd order two
embeddings, $X_t\hookrightarrow Y$ and $X_t \hookrightarrow
\left.E\right|_{X_t}$,
which are equivalent at 0th and 1st orders. We denote its dual
sheaf by~$\Delta_{X_t}^{[2]}(Y,E)$. The crucial property of this
creature is that there is a commutative diagram,
\begin{equation}
\begin{array}{rcccccl}
0 \ron & T_{t}M     & \ron & T_{t}^{[2]}M & \ron & S^2(T_{t}M) &\ron 0 \\
       &\downarrow&      &\downarrow  &      &\downarrow      &       \\
0\ron & H^0(X_t,N_t) &\ron &H^0\left(X_t,\Delta_{X_t}^{[2]}(Y,E)\right)&
\ron&H^0\left(X_t,S^2(N_t)\right)&\\
\end{array} \label{diagram}
\end{equation}
which relates the first floors of two towers of infinitesimal neighbourhoods
generated by embeddings $t\hookrightarrow M_0$ and $X_t\hookrightarrow Y$.
The verification is straightforward when one uses local coordinates
\cite{MeI}. Since, by assumption,
$H^1\left(X_t, N_t\otimes S^2( N^{*})\right)=0$, the exact sequence
\begin{equation}
0 \lon N_t \lon \Delta_{X_t}^{[2]}(Y,E) \lon S^2(N_t) \lon 0 \label{ext}
\end{equation}
splits. Moreover, its splitting is unique for the group
$H^0\left(X_t, N_t\otimes S^2(N_t^{*})\right)$ is also assumed to vanish.
The commutative diagram (\ref{diagram}) implies then a canonical splitting of
the exact sequence for $T^{[2]}_t M$ which is the same thing as a torsion-free
affine connection at $t$. This shows that the moduli space $M_0$ comes
equipped canonically with an induced torsion-free affine connection. The
fact that this connection is integrable on alpha surfaces is an easy
exercise. $\Box$
\vspace{2 mm}

For any given relative deformation problem~$X\hookrightarrow Y$, it is not
difficult to identify a large class of holomorphic vector bundles $E$ on the
ambient manifold $Y$ which have the property required by
Theorem~\ref{2.1}. Indeed, take any holomorphic distribution $E\subset TY$
on $Y$ which is transverse to $X$ and has rank equal to $\codim\, X$. Then, for
all $t$ in a sufficiently small neighbourhood $M_0$ of the point $t_0\in M$
corresponding to $X$, the submanifolds $X_t\hookrightarrow Y$ remain transverse
to $E$, and one has a {\em canonical} isomorphism $\left.E\right|_{X_t} =
N_{X_t}$.

Let us illustrate Theorem~\ref{2.1} by examples. The first one is just
the original Penrose non-linear graviton construction \cite{Penrose}.
Recall that Penrose
established a  one-to-one correspondence between ``small''
4-dimensional complex Riemannian manifolds $M_0$ with self-dual Riemann tensor
and 3-dimensional complex
manifolds $Y$ equipped with the following data:
(i) a submanifold $X\simeq \Bbb CP^{1}$ with normal bundle
$N\simeq\Bbb C^{2}\otimes\f_{\Bbb CP^{1}}(1)$;
(ii) a fibration $\pi: Y \lon X$; (iii) a global non-vanishing section
of the ``twisted'' determinant bundle $\det(V_{\pi})\otimes \pi^{*}(\f(2))$,
where $V_{\pi}$ is the vector bundle of $\pi$-vertical vector
fields. Let us look at this data from the point of view of Theorem~\ref{2.1}.
The distribution $V_{\pi}$ is clearly transversal
to $X$, and one easily checks that
$H^{0}\left(X_t, N_t\otimes S^{2}(N_t^{*})\right) =
H^{1}\left(X_t, N_t\otimes S^{2}(N_t^{*})\right) = 0$. Therefore, by
Theorem~\ref{2.1}, the data~(i) and~(ii) imply that there is
a torsion-free affine connection induced on $M_0$ which satisfies natural
integrability conditions. What is the role of datum (iii)? It will be shown
in section~6 that datum~(iii) ensures that this induced connection is
precisely the Levi-Civita connection of a holomorphic metric with
self-dual Riemann tensor. Therefore, the non-linear graviton construction is
one of the manifestations of the phenomenon envisaged by Theorem~\ref{2.1}.
Another manifestation is provided by relative deformation problems studied
in the context of quaternionic geometry. Let $Y$ be a complex
$(2n+1)$-dimensional manifold equipped with a holomorphic contact structure
which is a maximally non-degenerate rank $2n$ holomorphic distribution
$D\subset TY$. Let $X$ be a rational curve embedded into $Y$ transversely to
$D$ and  with normal bundle  $N={\Bbb C}^{2n}\otimes \f(1)$. Then
Theorem~\ref{2.1} says that the $4n$-dimensional Kodaira moduli space $M_0$
comes equipped with a unique torsion-free affine connection satisfying natural
integrability conditions. This is in accordance with well-known results in
twistor theory. The case $n=1$ has been studied by Ward~\cite{W2} and
Hitchin~\cite{Hitchin} who showed that $M_0$ has an induced complex
Riemannian metric satisfying self-dual Einstein equations with
non-zero scalar curvature. The case $n\geq 2$ has been investigated by
Bailey and Eastwood~\cite{BE}, LeBrun~\cite{suka}, and Pedersen and
Poon~\cite{PP}  who proved that the moduli space $M_0$ comes equipped
canonically with a torsion-free connection  compatible with
the induced (complexified) quaternionic structure on $M$.
\vspace{2 mm}

\paragraph{5. Families of affine connections on moduli space.}
What happens when $H^{0}\left(X_t, N_t\otimes S^{2}(N_t^{*})\right)\neq 0$?
In this case, as we noted in section~$3$, the locally free $\f_M$-module
$TM\otimes S^2(\OM)$ comes equipped with a non-zero $\f_M$-submodule $\blm$.
If $\nabla_1: TM\lon TM\otimes\OM$ and $\nabla_2: TM\lon TM\otimes\OM$ are
two torsion-free affine connections on $M$, then their difference,
$\nabla_1 - \nabla_2$, is a global section of $TM\otimes S^2(\OM)$. We
say that $\nabla_1$ and $\nabla_2$ are $\blm$-{\em equivalent}, if
$$
\nabla_1 -\nabla_2\in H^0(M,
\blm).
$$
We define a $\blm$-{\em connection} on $M$
as a collection of ordinary torsion-free affine connections
$\left\{\nabla_{i}\ | \ i\in I\right\}$ on an open covering $\left\{U_{i}\ | \
 i\in I\right\}$ of $M$ which, on overlaps $U_{ij} = U_{i} \cap U_{j}$, have
their differences in $H^{0}(U_{ij}, \blm)$. Locally, a $\blm$-connection
is the same thing as a $\blm$-equivalence
class of affine connections, but globally they are different ---
the obstruction for existence of a $\blm$-connection on $M$ lies
in $H^{1}\left(M,TM\otimes S^2(\OM)/\blm\right)$, while the obstruction
for existence of a $\blm$-equivalence
class of affine connections is an element
of $H^{1}\left(M,TM\otimes S^2(\OM)\right)$. A submanifold of~$M$
is said to be totally geodesic relative to a $\blm$-connection if it is
totally geodesic relative to
{\em each} of its local representatives $\nabla_{i}$.

\begin{theorem}\label{2.2}
Let \{$\X_{t}\hookrightarrow Y \mid t \in \M$\} be a complete family of compact
submanifolds.
Suppose there is a holomorphic vector bundle $E$ on $Y$ such that
$\left.E\right|_{X_t} = N_t$ for all $t$ in some domain $M_0\subseteq M$.
If $H^{1}\left(X_t, N_t\otimes S^{2}(N_t^{*})\right) = 0$ for each $t\in M_0$,
then $M_0$ comes equipped canonically with
an induced $\blm$-connection such that, for every $y\in Y'$, the associated
alpha subspace $\alpha_y\cap M_0$ is totally geodesic.
\end{theorem}
{\em Proof}. Since $H^1\left(X_t, N_t\otimes S^{2}(N_t^{*})\right)$, the exact
sequence (\ref{ext}) splits, and any such splitting, i.e.\ an
$\f_{X_t}$-morphism $\Delta_{X_t}^{[2]}(Y,E) \ron N_t \oplus S^2(N_t^{*})$,
induces an associated (via commutative diagram (\ref{diagram})) splitting
of the exact sequence for the second-order tangent space $T^{[2]}_t M$ which
in turn induces a torsion-free affine coonection at $t\in M_0$.
The set of all splittings of (\ref{ext}) is a principal homogeneous space
for the group $H^0\left(X_t, N_t\otimes S^{2}(N_t^{*})\right)$. Therefore,
the set of induced torsion-free affine connections at $t$ is a principal
homogeneous space for the group
$\blm_t$. $\Box$
\vspace{2 mm}

Next we study another equivalence class
of affine connections which is often induced on Kodaira moduli spaces.
Let $\sym$ denote the natural projection
$$
\sym:  TM\otimes\OM\otimes
\OM \rightarrow TM\otimes S^2(\OM) .
$$
Associated with a subsheaf $\Xi\subset TM\otimes\OM$
there is a concept of  $\Xi$-{\em connection}
on $M$ which is, by definition, a collection of ordinary torsion-free affine
connections $\left\{\nabla_{i}\ | \ i\in I\right\}$ on an open covering
$\left\{U_{i}\ | \ i\in I\right\}$ of $M$ which, on overlaps $U_{ij} =
U_{i} \cap U_{j}$, have
their differences in $H^{0}\left(U_{ij}, \sym(\Xi\otimes\OM)\right )$.
If $\Xi$ happens to be the structure sheaf $\f_M$
embedded diagonally into $TM\otimes\OM$, then a $\Xi$-connection is
nothing but a torsion-free projective connection on $M$.

\begin{theorem}\label{2.3}
Let \{$\X_{t}\hookrightarrow Y \mid t \in \M$\} be a complete family of compact
submanifolds.
Suppose there is a holomorphic vector bundle $E$ on $Y$ such that
$\left.E\right|_{X_t} = N_t$ for all $t$ in some domain $M_0\subseteq M$.
If $H^{1}(X_t, N_t\otimes N_t^{*}) = 0$ for each $t\in M_0$, then $M_0$ comes
equipped canonically with
an induced $\blm_{H^{0}(X,N\otimes N^{*})}$-connection such that, for every
$y\in Y'$, the associated
alpha subspace $\alpha_y\cap M_0$ is totally geodesic.
\end{theorem}
This theorem is proved by a slight modification of the construction
used in the proof of Theorem~1 \cite{MeI}. Note that the condition $H^{1}(X_t,
N_t\otimes N_t^{*}) = 0$ only says that
$N_t$ is a rigid vector bundle on $X_t$. There is an important case
when the requirement $\left.E\right|_{X_t} = N_t$ can be safely replaced by
a much weaker condition that vector bundles
$\left.E\right|_{X_t}$ and $N_t$ are simply equivalent to each other,
$\left.E\right|_{X_t} \simeq N_t$, which means that they define the same
cohomology class in $H^1(X_t, GL(p,\f_{X_t}))$, $p=\codim X$.

\begin{theorem}\label{2.4}
Let \{$\X_{t}\hookrightarrow Y \mid t \in \M$\} be a complete family of
compact 1-codimensional submanifolds.
Suppose there is a holomorphic line bundle $E$ on $Y$ such that
$\left.E\right|_{X_t} \simeq N_t$ for all $t$ in some domain $M_0\subseteq M$.
If $H^{1}(X_t, \f_{X_t}) = 0$ for each $t\in M_0$, then $M_0$ comes
equipped canonically with
an induced torsion-free projective connection such that, for every
$y\in Y'$, the associated
alpha subspace $\alpha_y\cap M_0$ is totally geodesic.
\end{theorem}
{\em Proof}. Fix any isomorphism $\phi: \left.E\right|_{X_t} \lon N_t$
and apply  Theorem~\ref{2.4} to construct a projective connection $\nabla$
on $M_0$. Due to compactness of $X_t$, the isomorphism $\phi$ is defined up to
multiplication by a non-zero constant. The point is that the construction of
$\nabla$  scetched in the proof of Theorem~\ref{2.1} is invariant under
such a transformation which shows that $\nabla$ is actually independent
on a particular choice of the isomorphism $\phi$ used in
its construction.
$\Box$

\begin{corollary}\label{2.5}
Let \{$\X_{t}\hookrightarrow Y \mid t \in \M$\} be a complete family of
compact 1-codimensional submanifolds. If $H^{1}(X_{t_0}, \f_{X_{t_0}}) = 0$
for some  $t_0\in M$, then there is an open neigbourhood $M_0$ of $t_0$
which comes
equipped canonically with
an induced torsion-free projective connection such that, for every
$y\in Y'$, the associated
alpha subspace $\alpha_y\cap M_0$ is totally geodesic.
\end{corollary}
{\em Proof}. By semi-continuity principle, there is an open neighbourhood
$U$ of $t_0$ such that $H^1(X_t, \f_{X_t})=0$ for each $t\in U$. If
$L$ is the divisor line bunlde of $X_{t_0}$, then
$\left.L^{*}\right|_{X_t} \simeq N_t$ for all $t$
in probably smaller neighbourhood $M_0$. This fact combined with
Theorem~\ref{2.4} implies the desired result. $\Box$

An independent (and more effective from the viewpoint of computing
induced projective structures) proof of this Corollary is given in \cite{MP}.
\vspace{3 mm}

\paragraph{6. Fibred complex manifolds.}
Suppose that the ambient manifold $Y$ has the structure of a holomorphic
fibration over its compact submanifold $X$, i.e.\ there is a submersive
holomorphic map $\pi: Y\lon X$.
If $H^1(X,N)=0$, then, by Kodaira's theorem, $X$ belongs to the complete
family $\{X_t \hookrightarrow Y \mid t\in M\}$ of compact submanifolds.
The submanifold
$X$ is transverse to the distribution $V_{\pi}$ of $\pi$-vertical vector fields
on $Y$, and so is $X_t$ for every $t$ in some neighbourhood $M_0$ of the point
$t_0\in M$ which corresponds to $X$.
 Thus $Y$ has the structure of a holomorphic fibration $Y\lon X_t$
for $t\in M_0$. If also $H^1\left(X_t, N_t\otimes S^{2}(N_t^{*})\right)=0$,
then, by Theorem~\ref{2.2}, the moduli space comes equipped with a family
of induced torsion-free affine connections. Holonomy groups of these
induced connections can be estimated with surprising ease.
\begin{theorem}\label{hol 1}
Let $\nabla$ be an induced connection on the moduli space $M_0$.
If the function
$$
f: t\lon \dim H^0(X_t,N_t\otimes N_t^{*})
$$
is constant on $M_0$, then the holonomy algebra of
$\nabla$ is a subalgebra of the finite dimensional Lie algebra
{\em $H^{0}(X, N\otimes N^{*})$}.
If in addition there is a holomorphic line bundle $L$ on $X$
such that the bundle {\em $\pi^{*}(L)\otimes \det\,V_{\pi}^{*}$} admits a
nowhere vanishing holomorphic section, then the holonomy algebra of $\nabla$
is a subalgebra of $H^0(X, N\otimes_0 N^{*})$, where $\otimes_0$ denotes
trace-free tensor product.
\end{theorem}
This theorem follows from the correspondence between {\em third}\, floors
of the two towers of infinitesimal neighbourhoods of embeddings
$t\hookrightarrow M$ and $X_t\hookrightarrow Y$
which has been studied in \cite{Me1}.

{\em Examples}. 1. Let $X= \Bbb C\Bbb P^1$ be the
projective line embedded into a
3-dimensional complex manifold $Y$ with normal bundle $N\simeq
\Bbb C^2\otimes \f (1)$. If $Y$ has the structure of a holomorphic fibration
over $X$, then Theorem~\ref{2.1} says that there is
an induced connection $\nabla$ on the moduli space $M_0$. By
Theorem~\ref{hol 1}, the holonomy algebra of $\nabla$ is contained in $H^0(X,
N\otimes
N^{*}) = M_{2,2}(\Bbb C)\subset co(4,\Bbb C)$, where $M_{2,2}(\Bbb C)$ is the
space of $2\times 2$ complex
matrices and $co(4,\Bbb C)= sl(2,\Bbb C) + sl(2,\Bbb C) + \Bbb C$ the
complexified conformal algebra. Since $\nabla$ is torsion-free, this fact
implies that $\nabla$
is a complex Weyl connection on the 4-dimensional complex conformal manifold
$M_0$ which
has the anti-self-dual
parts of the Weyl tensor and the antisymmetrized Ricci tensor vanishing and
satisfies the Einstein-Weyl equations. In fact any such  connection arises
locally
in this way \cite{Penrose,Boyer,PS}.

2. Let the pair $\left(X= \Bbb C\Bbb P^1, Y\right)$ be the same
as in Example~1 and assume that the bundle $\det V_{\pi}^{*}\otimes
 \pi^{*}\left(\f (2)\right)$  admits a nowhere
vanishing holomorphic section. By Theorem~\ref{hol 1}, the holonomy group
of the induced connection $\nabla$ on the 4-dimensional parameter space $M_0$
is a subgroup of $SL(2,\Bbb C)\subset SO(4,\Bbb C)=SL(2,\Bbb C)\times_{\Bbb
Z_2} SL(2,\Bbb C)$.
This means that $\nabla$ is the Levi-Civita
connection of a complex Riemannian metric on $M_0$ which is Ricci-flat
and has the anti-self-dual part of the Weyl tensor vanishing. Again any such
connection arises locally in this way \cite{Penrose}.

3. Let $X= \Bbb C\Bbb P^1$ be the projective line embedded into a
$(2k+1)$-dimensional complex manifold $Y$ with normal bundle $N\simeq
\Bbb C^{2k}\otimes \f (1)$, $k\geq 2$. The Kodaira moduli space $M$ is then
a $4k$-dimensional complex manifold possessing a complexified almost
quaternionic structure. If $Y$ has the structure of a holomorphic fibration
over $X$,
then, by Theorems~\ref{2.1} and~\ref{hol 1}, there is
an induced  connection $\nabla$ on $M_0\subset
M$ with holonomy in $GL(2k,\Bbb C)$
which implies that $\nabla$ is a complexified Obata connection. It is not
difficult
to show using results of Bailey and Eastwood \cite{BE} that any local
complexified Obata connection can be constructed in this way. A twistor
characterization of real analytic Obata connections
can be obtained from this holomorphic picture in the usual way
(cf.\  \cite{Manin}).

\end{document}